\documentclass{amsart}
\usepackage{amsmath}%
\usepackage{amsfonts}%
\usepackage{amssymb}%
\usepackage{graphicx}
\theoremstyle{plain}

\numberwithin{equation}{section}
\begin{document}
\title[degrees for eigensystems of rational self-maps]{The dynamical and arithmetical degrees for eigensystems of rational self-maps}
\author{Jorge Mello}
\address[]
{Universidade Federal do Rio de Janeiro, Instituto de Matem\'{a}tica. mailing adress:\newline Rua Aurora 57/101, Penha Circular, 21020-380 Rio de Janeiro, RJ, Brasil. %
 } \email[]{jmelloguitar@gmail.com}%
\urladdr{https://sites.google.com/site/algebraufrj/}

\thanks{The author would like to say thanks to J. Silverman and S. Kawaguchi.}
\date{december 20, 2017}
\dedicatory{}

\begin{abstract}
We define arithmetical and dynamical degrees for dynamical systems with several rational maps on projective varieties, study their properties and relations, and prove the existence of a canonical height function associated with divisorial relations in the N\'{e}ron-Severi Group over Global fields of characteristic zero, when the rational maps are morphisms. For such, we show that for any Weil height $h_X$ with respect to an ample divisor on a projective variety $X$, any dynamical system $\mathcal{F}$ of rational self-maps on $X$, and any $\epsilon>0$, there is a positive constant $C=C(X, h_X, f, \epsilon)$ such that \begin{center}
$\sum_{f \in \mathcal{F}_n} h^+_X(f(P)) \leq C. k^n.(\delta_{\mathcal{F}} + \epsilon)^n .  h^+_X(P)$
\end{center} for all points $P$ whose $\mathcal{F}$-orbit is well defined.
\end{abstract}
\maketitle

\section{Introduction}
Weil heights have an important role in Diophantine geometry, and particular Weil heights with nice properties, called canonical heights, are sometimes very useful. The theory of canonical heights has had deep applications throughout the field of Arithmetic geometry.

Over abelian varieties $A$ defined over a number field $K$, N\'{e}ron and Tate constructed canonical height functions $\hat{h}_L: A(\bar{K}) \rightarrow \mathbb{R}$ with respect to symmetric ample line bundles $L$ which enjoy nice properties, and can be used to prove Mordell-Weil theorem for the rational points of the variety. More generally, in [4], Call and Silverman constructed canonical height functions on projective varieties $X$ defined over a number field which admit a morphism $f:X \rightarrow X$ with $f^*(L) \cong L^{\otimes d}$ for some line bundle $L$ and some $d >1$. In another direction, Silverman [19] constructed canonical height functions on certain $K3$ surfaces $S$ with two involutions $\sigma_1, \sigma_2$ (called Wheler's $K3$ surfaces) and developed an arithmetic theory analogous to the arithmetic theory on abelian varieties.

It was an idea of Kawaguchi [10] to consider polarized dynamical systems of several maps, namely, given $X/K$ a projective variety, $f_1,...f_k:X \rightarrow X$ morphisms on defined over $K$, $\mathcal{L}$ an invertible sheaf on $X$ and a real number $d>k$ so that $f_1^*\mathcal{L}\otimes ... \otimes f_k^*\mathcal{L} \cong \mathcal{L}^{\otimes d}$, he constructed a canonical height function associated to the polarized dynamical system $(X, f_1,..., f_k, \mathcal{L})$ that generalizes the earlier constructions mentioned above. In the Wheler's K3 surfaces' case above, for example, the canonical height defined by Silverman arises from the system formed by $(\sigma_1, \sigma_2)$ by Kawagushi's method.

Given $X/\mathbb{C}$ smooth projective variety, $f:X \dashrightarrow X$ dominant rational map inducing $f^*:$NS$(X)_{\mathbb{R}} \rightarrow$NS$(X)_{\mathbb{R}}$ on the N\'{e}ron-Severi group, the dynamical degree is defined as $\delta_f := \lim_{n \rightarrow \infty} \rho((f^n)^*)^{\frac{1}{n}}$, where $\rho$ denotes the spectral radius of a given linear map, or the biggest number among the absolute values of its eigenvalues. This limit converges and is a birational invariant that has been much studied over the last decades. In [12] we find a list of references.

In [12], Kawaguchi and Silverman studied an analogous arithmetic degree for $X$ and $f$ defined over $\bar{\mathbb{Q}}$ on points with well defined foward orbit over $\bar{\mathbb{Q}}$. Namely, $\alpha_f(P):= \lim_{n \rightarrow \infty} h^+_X(f^n(P))^{\frac{1}{n}}$, where $h_X$ is a Weil height relative to an ample divisor and $h^+_X= \max \{1, h_X\}$. Such degree measures the arithmetic complexity of the orbit of $P$ by $f$, and $\log \alpha_f(P)$ has been interpreted as a measure of the arithmetic entropy of the orbit $\mathcal{O}_f(P)$. It is showed in [12] that the arithmetic degree determines the height counting function for points in orbits, and that the arithmetic complexity of the $f$-orbit of an algebraic point never exceeds the geometrical-dynamical complexity of the map $f$, as well as more arithmetic consequences. We ask if this kind of research could be done in the setting of general dynamical systems as treated by Kawaguchi, with several maps, as in the case of Wheler's K3 surfaces. This is the first subject found in this work.

Given $X/K$ be a projective variety, $f_1,...,f_k:X \dashrightarrow X$ rational maps, $\mathcal{F}_n=\{f_{i_1} \circ ... \circ f_{i_n} ; i_j =1,...,k \}$, we define a more general dynamical degree of a system of maps as $\delta_{\mathcal{F}}=\lim \sup_{n\rightarrow \infty}\max_{f \in \mathcal{F}_n} \rho (f^*)^{\frac{1}{n}}$, and extend the definition of arithmetic degree for $\alpha_{\mathcal{F}}(P)= \frac{1}{k} \lim_{n \rightarrow \infty}\{ \sum_{f \in \mathcal{F}_n} h_X^+(f(P))\}^{\frac{1}{n}}$, obtaining also the convergence of $\delta_{\mathcal{F}}$, and that $\alpha_{\mathcal{F}}(P) \leq \delta_{\mathcal{F}}$ when $\alpha_{\mathcal{F}}(P)$ exists. Motivated by [12], we give an elementary proof that our new arithmetic degree is related with height counting functions in orbits, when $\alpha_{\mathcal{F}}(P)$ exists, by:
\begin{center}
$\lim_{B \rightarrow \infty} \dfrac{ \# \{ n \geq 0 ; \sum_{f \in \mathcal{F}_n} h_X(f(P)) \leq B\}}{ \log B}= \dfrac{1}{ \log (k. \alpha_{\mathcal{F}}(P))} $,
\end{center}
\begin{center}
$\lim \inf _{B \rightarrow \infty} (\# \{Q \in \mathcal{O}_{\mathcal{F}}(P); h_X(Q) \leq B \})^{\frac{1}{\log B}} \geq  k^{ \frac{1}{ \log (k. \alpha_{\mathcal{F}}(P)})}$.
\end{center}
We are able to extend theorem 1 of [12], showing explicitely how the dynamical degree of a system with several maps can offer an uniform upper bound for heights on iterates of points in orbits, when $K$ is a number field or an one variable function field. Precisely, for every $\epsilon >0$,
there exists a positive constant $C=C(X, h_X, f, \epsilon)$ such that for all $P \in X_{\mathcal{F}}(\bar{K})$ and all $n \geq 0$,
\begin{center}
$\sum_{f \in \mathcal{F}_n} h^+_X(f(P)) \leq C. k^n.(\delta_{\mathcal{F}} + \epsilon)^n .  h^+_X(P).$
\end{center} In particular, $h^+_X(f(P)) \leq  C. k^n.(\delta_{\mathcal{F}} + \epsilon)^n .  h^+_X(P)$ for all $f \in \mathcal{F}_n.$

This theorem becomes a tool to show the second very important theorem of this work. As we have seen, for a pair $(X/K, f_1,...,f_k, L)$ with $k$ self-morphisms on $X$ over $K$, and $L$ a divisor satisfying a linear equivalence $\otimes^k_{i=1}f^*_i(L) \sim L^{\otimes d}$ for $d>k$, there is a well known theory of canonical heights developed by Kawaguchi in [10]. Now we are partially able to generalize this to cover the case that the relation $\otimes^k_{i=1}f^*_i(L) \equiv L^{\otimes d}$ is only an algebraic relation. Hence the limit 
\begin{center}
$\hat{h}_{L,\mathcal{F}}(P)= \lim_{n \rightarrow \infty}\dfrac{1}{d^n}\sum_{ f \in \mathcal{F}_n} h_L(f(P)).$
\end{center}
converges for certain eigendivisor classes relative to algebraic relation. For $L$ ample and $K$ a number field, we obtain that :
\begin{center}
$\hat{h}_{L,\mathcal{F}}(P)=0 \iff  P$ has finite $\mathcal{F}$-orbit.
\end{center}
These kind of generalization was firstly done for just one morphism by Y. Matsuzawa in [15], extending Call and Silverman's theory of canonical heights in [4], and we work out for several maps in the present work.

\section{Notation, and first definitions}

Throughout this work, $K$ will be either a number field or a one-dimensional function field of characteristic 0 . We let $\bar{K}$ be an algebraic closure of $K$. The uple $ (X, f_1,...,f_k)$ is called a dynamical system, where either $X$ is a smooth projective variety and $f_i:X \dashrightarrow X$ are dominant rational maps all defined over $K,$  or $X$ is a normal projective variety and $f_i:X \dashrightarrow X$ are dominant morphisms.

We denote  by $ h_X:X(\bar{K})\rightarrow [0, \infty)$ the absolute logarithmic Weil height function relative to an ample divisor $A$ of $X$, and for convenience we set $h_X^+(P)$ to be $ \max \{1, h_X(P)\}.$

The sets of iterates of the maps in the system are denoted by $\mathcal{F}_0=\{ $Id$\},  \mathcal{F}_1= \mathcal{F} =\{f_1,...,f_k\}$, and $\mathcal{F}_n=\{f_{i_1} \circ ... \circ f_{i_n} ; i_j =1,...,k \}$, inducing what we call $\mathcal{O}_{\mathcal{F}}(P)$ the forward $\mathcal{F}$-orbit of $P$=$\{ f(P); f \in \bigcup_{n \in \mathbb{N}} \mathcal{F}_n \}.$ A point $P$ is said preperiodic when its $\mathcal{F}$-orbit is a finite set.

We write $ I_{f_i}$ for the indeterminacy locus of $f_i$, i.e., the set of points which $f_i$ is not well-defined, and $ I_{\mathcal{F}}$ for $ \bigcup_{i=1}^k  I_{f_i}$. Also we define $ X_{\mathcal{F}}(\bar{K})$ as the set of points $P \in X(\bar{K})$ whose forward orbit is well-defined, in other words, $\mathcal{O}_{\mathcal{F}}(P) \cap  I_{\mathcal{F}} = \emptyset$. 

The set of Cartier divisors on $X$ is denoted by Div$(X)$, while Pic$(X)$ denotes The Picard group of $X$, and NS$(X)=\mbox{Pic}(X)/\mbox{Pic}^{0}(X)$ is called the Neron-Severi Group of $X$. The equality in this group is denoted by the symbol $\equiv$, which is called algebraic equivalence.

Given a rational map $f:X \dashrightarrow X$, the linear map induced on the tensorized N\'{e}ron-Severi Group NS$(X)_{\mathbb{R}}=$ NS$(X) \otimes \mathbb{R}$ is denoted by $f^*$. So, when looking for a dynamical system $(X,\mathcal{F})$, it is convenient for us to use the notation  $\rho(\mathcal{F}_n):= \max_{f \in \mathcal{F}_n} \rho (f^*,$NS$(X)_{\mathbb{R}})$.

For definitions and properties about Weil height functions, we refer to [8].

Next, we define the dynamical degree of a set of rational maps on a complex variety, which is a measure of the geometric complexity of the iterates of the maps in the set, when it exists. This is a generalization for several morphisms of the dynamical degree appearing in of [12]. \newline 

{\bf Definition 2.1: } 
\textit{Let $X/ \mathbb{C}$ be a (smooth) projective variety and let $\mathcal{F}$ be as above. The dynamical degree of $\mathcal{F}$, when it exists, is defined by}
\begin{center}
$\delta_{\mathcal{F}}=\lim \sup_{n\rightarrow \infty}\rho(\mathcal{F}_n)^{\frac{1}{n}}$
\end{center} 

In this sense, we also generalize the second definition in the introduction of [12], introducing now the arithmetic degree of a system of maps $\mathcal{F}$ at a point $P$. This degree measures the growth rate of the heights of $n$-iterates of the point by maps of the system as $n$ grows, and so it is a measure of the arithmetic complexity of $\mathcal{O}_{\mathcal{F}}(P)$.\newline

{\bf Definition 2.2: } 
\textit{Let $P \in  X_{\mathcal{F}}(\bar{K}).$ The arithmetic degree of $\mathcal{F}$ at $P$ is the quantity}
\begin{center}
$\alpha_{\mathcal{F}}(P)= \frac{1}{k} \lim_{n \rightarrow \infty}\{ \sum_{f \in \mathcal{F}_n} h_X^+(f(P))\}^{\frac{1}{n}} $
\end{center} \textit{assuming that the limit exists.}\newline

{\bf Definition 2.3: } \textit{In the lack of the convergence, we define the upper and  the lower arithmetic degrees as}
\begin{center}
$\bar{\alpha}_{\mathcal{F}}(P)= \frac{1}{k} \lim \sup_{n \rightarrow \infty} \{ \sum_{f \in \mathcal{F}_n} h_X^+(f(P))\}^{\frac{1}{n}} $ 

$\underline{\alpha}_{\mathcal{F}}(P)= \frac{1}{k} \lim \inf_{n \rightarrow \infty} \{ \sum_{f \in \mathcal{F}_n} h_X^+(f(P))\}^{\frac{1}{n}} $
\end{center}

{\bf Remark 2.4: } 
Let $X$ be a projective variety and $D$ a Cartier divisor. If \newline  $f:X \rightarrow X $ is a surjective morphism, then $f^*D$ is a Cartier divisor. In the case where $X$ is smooth, and $f:X \dashrightarrow X$ a merely rational map, we take a smooth projective variety $\tilde{X}$ and a birational morphism $\pi: \tilde{X} \rightarrow X$ such that $\tilde{f} := f \circ \pi: \tilde{X} \rightarrow X$ is a morphism. And we define $f^*D:=\pi_*(\tilde{f}^*D).$ It is not hard to verify that this definition is independent of the choice of $X$ and $\pi$. This is done in section 1 of [12] for example.

\section{First properties for the arithmetic degree}

In this section we check that the upper and lower degrees defined in the end of the section above are independent of the Weil height function chosen for $X$, and so they are well defined. Some examples of these degrees are computed is this section as well. We also present and prove our first counting result for points in orbits for several maps, and state an elementary and useful linear algebra's lemma.\newline

{\bf Proposition 3.1: }
\textit{The upper and lower arithmetic degrees $\bar{\alpha}_{\mathcal{F}}(P)$ and $ \underline{\alpha}_{\mathcal{F}}(P)$ are independent of the choice of the height function $h_X$.}

\begin{proof}
If the $\mathcal{F}$-orbit of $P$ is finite, then the limit $\alpha_{\mathcal{F}}(P)$ exists and is equal to 1, by definition of such limit, whatever the choice of $h_X$ is. So we consider the case when $P$ is not preperiodic, which allows us to replace $h_X^+$ with $h_X$ when taking limits. 

Let $h$ and $h^{\prime}$ be the heights induced on $X$ by ample divisors $D$ and $D^{\prime}$ respectively, and let the respective arithmetic degrees denoted by $\bar{\alpha}_{\mathcal{F}}(P)$,  $\underline{\alpha}_{\mathcal{F}}(P)$,  ${\bar{\alpha}}^{\prime}_{\mathcal{F}}(P)$, ${\underline{\alpha}}^{\prime}_{\mathcal{F}}(P)$. By the definition of ampleness, there is an integer $m$ such that $mD- D^{\prime}$ is ample, and thus the functorial properties of height functions imply the existence of a non-negative constant $C$ such that:
\begin{center}
$mh(Q) \geq h^{\prime}(Q) - C $ for all $ Q \in X(\bar{K}).$ 
\end{center}
We can choose a sequence of indices $\mathcal{N} \subset \mathbb{N}$ such that:
\begin{center}
$\bar{\alpha}^{\prime}_{\mathcal{F}}(P)= \frac{1}{k} \lim \sup_{n \rightarrow \infty} \{ \sum_{f \in \mathcal{F}_n} h^{\prime}(f(P))\}^{\frac{1}{n}}  = \frac{1}{k} \lim_{n \in \mathcal{N}} \{ \sum_{f \in \mathcal{F}_n} h^{\prime}(f(P))\}^{\frac{1}{n}}$ 
\end{center}
Then \newline \newline
$\bar{\alpha}^{\prime}_{\mathcal{F}}(P)=\frac{1}{k} \lim_{n \in \mathcal{N}} \{ \sum_{f \in \mathcal{F}_n} h^{\prime}(f(P))\}^{\frac{1}{n}}  \newline  \newline  \leq \frac{1}{k} \lim_{n \in \mathcal{N}} \{ \sum_{f \in \mathcal{F}_n}m h(f(P))+C \}^{\frac{1}{n}}  ~~\newline \newline \leq \frac{1}{k} \lim \sup_{n \rightarrow \infty} \{ \sum_{f \in \mathcal{F}_n}m h(f(P))+C \}^{\frac{1}{n}} ~  \newline \newline
=  \frac{1}{k} \lim \sup_{n \rightarrow \infty} \{ m(\sum_{f \in \mathcal{F}_n} h(f(P)))+Ck^n \}^{\frac{1}{n}} \newline \newline
= \frac{1}{k} \lim \sup_{n \rightarrow \infty} \{ \sum_{f \in \mathcal{F}_n} h(f(P))\}^{\frac{1}{n}} \newline \newline =
\bar{\alpha}_{\mathcal{F}}(P). $ \newline

This proves the inequality for the upper arithmetic degrees. Reversing the roles of $h$ and $h^{\prime}$ in the calculation above we also prove the opposite inequality, which demonstrates that  $\bar{\alpha}_{\mathcal{F}}(P)={\bar{\alpha}}^{\prime}_{\mathcal{F}}(P).$ In the same way we prove that  $\underline{\alpha}_{\mathcal{F}}(P)={\underline{\alpha}}^{\prime}_{\mathcal{F}}(P).$
\end{proof}
Our next lemma says that points belonging to a fixed orbit have their upper and lower arithmetic degrees    bounded from above by the respective arithmetic degrees of the given orbit generator point. \newline \newline
{\bf Lemma 3.2: }
\textit{Let $\mathcal{F}=\{f_1,..., f_k\}$ be a set of self-rational maps on $X$ defined over $\bar{K}$. Then, for all $P \in X_{\mathcal{F}}(\bar{K})$, all $l \geq 0$, and all $g \in \mathcal{F}_l$, }
\begin{center}
$\bar{\alpha}_{\mathcal{F}}(g(P)) \leq \bar{\alpha}_{\mathcal{F}}(P)$ and  ~$\underline{\alpha}_{\mathcal{F}}(g(P)) \leq \underline{\alpha}_{\mathcal{F}}(P)$
\end{center}
\begin{proof} We calculate \newline \newline
$\bar{\alpha}_{\mathcal{F}}(g(P)) = \frac{1}{k} \lim \sup_{n \rightarrow \infty} \{ \sum_{f \in \mathcal{F}_n} h_X^+(f(g(P)))\}^{\frac{1}{n}}\newline \newline
  \frac{1}{k} \lim \sup_{n \rightarrow \infty} \{ \sum_{f \in \mathcal{F}_n , g^{\prime} \in \mathcal{F}_l } h_X^+(f(g^{\prime}(P))) - \sum_{f \in \mathcal{F}_n , g^{\prime} \in \mathcal{F}_l -\{g\}} h_X^+(f(g^{\prime}(P)))\}^{\frac{1}{n}} ~~~~~~~~~~~~~~~~~~~\newline  \newline
\leq \frac{1}{k} \lim \sup_{n \rightarrow \infty} \{ [\sum_{f \in \mathcal{F}_{n+l}} h_X^+(f(P))]+ O(1).k^{n+l}\}^{\frac{1}{n}}~~~~~~~~~~~~~~~~~~~~~~~~~~~~~~~~~~~~~~~~~~~~~~~~~~\newline \newline
=\frac{1}{k} \lim \sup_{n \rightarrow \infty} \{ \sum_{f \in \mathcal{F}_{n+l}} h_X^+(f(P))\}^{\frac{1}{n}}~~~~~~~~~~~~~~~~~~~~~~~~~~~~~~~~~~~~~~~~~~~~~~~~~~~~~~~~~~~~~~~~~~\newline \newline
=\frac{1}{k} \lim \sup_{n \rightarrow \infty} \{ \sum_{f \in \mathcal{F}_{n+l}} h_X^+(f(P))\}^{\frac{1}{n+l} . (1+ \frac{l}{n})} ~~~~~~~~~~~~~~~~~~~~~~~~~~~~~~~~~~~~~~~~~~~~~~~~~~~~~~~~~~~~~~~\newline \newline
=\frac{1}{k} \lim \sup_{n \rightarrow \infty} \{ \sum_{f \in \mathcal{F}_{n+l}} h_X^+(f(P))\}^{\frac{1}{n+l}} ~~~~~~~~~~~~~~~~~~~~~~~~~~~~~~~~~~~~~~~~~~~~~~~~~~~~~~~~~~~~~~~~~~~~\newline \newline
=\bar{\alpha}_{\mathcal{F}}(P)$~~~~~~~~~~~~~~~~~~~~~~~~~~~~~~~~~~~~~~~~~~~~~~~~~~~~~~~~~~~~~~~~~~~~~~~~~~~~~~~~~~~~~~~~~~~~~~~~~~~~~~~~~~~~~~~~~~~~~~~~~~~~~~~~~~~~~~~~~~~~~~~~~~~~~~~
\newline \newline  The proof for $\underline{\alpha}_{\mathcal{F}}(P)$ is similar.

\end{proof}

\newpage Here are some examples: \newline
 
$\bf{Example \: 3.3}$: Let $S$ be a $K3$ surface in $\mathbb{P}^2 \times \mathbb{P}^2$ given by the intersection of two hypersurfaces of bidegrees (1,1) and (2,2) over $\overline{\mathbb{Q}}$, and assume that NS$(S)  \cong \mathbb{Z}^2,$ generated by $L_i:= p_i^*O_{\mathbb{P}^2}(1), i=1,2$, where $p_i : S \rightarrow \mathbb{P}^2$ is the projection to the $i$-factor for $i=1,2.$ These induce  noncommuting involutions $\sigma_1, \sigma_2  \in $ Aut$(S)$. By [19, Lemma 2.1], we have \begin{center} $\sigma_i^*L_i \cong L_i, \sigma_i^*L_j \cong 4L_i - L_j,$ for $ i \neq j.$\end{center} The line bundle $L:= L_1 + L_2$ is ample on $S$ and satisfies $\sigma_1^*L + \sigma^*_2L \cong 4L$, and thus $h:= \hat{h}_{L, \{ \sigma_1, \sigma_2\}}$ exists on $S(\overline{\mathbb{Q}})$ by [10, theorem 1.2.1]. Noting that 
  $$ \sigma_1^* \sim \begin{bmatrix}
  1 & 4 \\
  0 & -1
 \end{bmatrix}, 
 \sigma_2^* \sim \begin{bmatrix}
  -1 & 0 \\
  4 & 1
 \end{bmatrix}, 
 (\sigma_1 \circ \sigma_2)^* \sim \begin{bmatrix}
  -1 & -4 \\
  4 & 15
 \end{bmatrix}, 
 (\sigma_2 \circ \sigma_1)^* \sim \begin{bmatrix}
  15 & 4 \\
  -4 & -1
 \end{bmatrix} , $$
$$
 (\sigma_1 \circ \sigma_2 \circ \sigma_1)^* \sim \begin{bmatrix}
  15 & 56 \\
  -4 & -15
 \end{bmatrix} , 
 (\sigma_2 \circ \sigma_1 \circ \sigma_2)^* \sim \begin{bmatrix}
  -15 & -4 \\
  56 & 15
 \end{bmatrix}, 
$$  we calculate that \begin{center}$\rho(\sigma_1^*)= 2+ \sqrt{3}, \rho( \sigma_2^*)= 2 + \sqrt{3}, \rho((\sigma_1 \circ \sigma_2 )^*)= 7 + 4 \sqrt{3}, \newline \rho(( \sigma_2 \circ \sigma_1)^*) = 7 + 4 \sqrt{3}, \rho((\sigma_1 \circ \sigma_2 \circ \sigma_1)^*)= 1, \rho((\sigma_2 \circ \sigma_1 \circ \sigma_2)^*)=1.$ \end{center} This gives that $\delta_{\{\sigma_1, \sigma_2\}}= 2 + \sqrt{3}. $ Furthermore, since $h$ is a Weil Height with respect to an ample divisor, 

\begin{center}$\alpha_{\{\sigma_1, \sigma_2\}}(P) = (1/2) . \lim_{n \rightarrow \infty } [\sum_{f \in \{ \sigma_1, \sigma_2 \}_n} h(f(P))]^{\frac{1}{n}}=1/2.[4^n.h(P)]^{ \frac{1}{n}}=2 $ \end{center} for all $P \in S(\bar{\mathbb{Q}})$ non-preperiodic, i.e, $P$ such that $h(P) \neq 0.$

Observe that in this case $\bar{\alpha}_{\{\sigma_1, \sigma_2\}}(P)=2 \leq 2 + \sqrt{3}= \delta_{\{\sigma_1, \sigma_2\}},$ which we will prove in Corollary 1.16 to be true in our general conditions. \newline

$\bf{Example \: 3.4}$: Let $S$ be a $K3$ surface in $\mathbb{P}^2 \times \mathbb{P}^2$, as in the example 1.4.5 of [10], given by the intersection of two hypersurfaces of bidegrees (1,2) and (2,1) over $\overline{\mathbb{Q}}$, and assume that NS$(S)  \cong \mathbb{Z}^2,$ generated by $L_i:= p_i^*O_{\mathbb{P}^2}(1), i=1,2$, where $p_i : S \rightarrow \mathbb{P}^2$ is the projection to the $i$-factor for $i=1,2.$ These induce  noncommuting involutions $\sigma_1, \sigma_2  \in $ Aut$(S)$. By similar computations we have $\sigma_i^*L_i \cong L_i, \sigma_i^*L_j \cong 5L_i - L_j,$ for $ i \neq j.$ The ample line bundle $L:= L_1 + L_2$ exists on $S$ and satisfies $\sigma_1^*L + \sigma^*_2L \cong 5L$, and thus $h:= \hat{h}_{L, \{ \sigma_1, \sigma_2\}}$ exists on $S(\overline{\mathbb{Q}})$ by [10, theorem 1.2.1]. Proceeding in the same way as in the previous example, we have that  \begin{center} $\bar{\alpha}_{\{\sigma_1, \sigma_2\}}(P) =5/2  \leq \sqrt{\dfrac{23 + 5 \sqrt{21}}{2}}={\delta}_{\{\sigma_1, \sigma_2\}}.$ \end{center}

$\bf{Example \: 3.5}$: Let $S$ be a hypersurface of tridegree (2,2,2) in $\mathbb{P}^1 \times \mathbb{P}^1 \times \mathbb{P}^1$ over $\overline{\mathbb{Q}}$, as in the example 1.4.6 of [10]. For $i=1,2,3,$ let $p_i:S \rightarrow \mathbb{P}^1 \times \mathbb{P}^1$ be the projection to the $(j,k)-$th factor with $\{i,j,k\}= \{1,2,3\}.$ Since $p_i$ is a double cover, it gives an involution $\sigma_i \in $ Aut$(S).$ Let also, $q_i:S \rightarrow \mathbb{P}^1$ be the projection to the $i-$th factor, and set $L_i := q_i^* O_{\mathbb{P}^1}, L:= L_1 + L_2 + L_3$ ample, and we assume that NS$(S) =  <L_1,L_2,L_3> \cong \mathbb{Z}^3.$ By similar computations as above we have
\begin{center}
$ \sigma_i^*(L_i) \cong -L_i +2L_j + 2 L_k $ for $ \{i,j,k\}= \{ 1,2,3\}\newline
\sigma_j^*(L_i) \cong L_i$ for $i \neq j.$
\end{center} Then $\sigma_1^*L + \sigma_2^*L + \sigma_3^*L \cong 5L,$ which gives us the existence of $h := \hat{h}_{L, \{ \sigma_1, \sigma_2, \sigma_3\}}$ by [10, theorem 1.2.1]. We note that if $h(P) \neq 0$, then a similar computation as in the previous examples yields $\alpha_{\{ \sigma_1,\sigma_2,\sigma_3\}}(P)= 5/3 $. While we can also calculate that: 
$$
 (\sigma_3 \circ \sigma_2 \circ \sigma_1)^* \sim \begin{bmatrix}
    1 & -2 & -2 \\
    2 & 3 & 10 \\
    2 & 6 & 15
  \end{bmatrix}
$$ with its big eigenvalue being aproximatelly $\rho( (\sigma_3 \circ \sigma_2 \circ \sigma_1)^*) \sim 18,3808$. As $(18,3808)^{1/3} \sim 2,639$, we have that $\delta_{\{\sigma_1,\sigma_2,\sigma_3\}} \geq  2,63 > 5/3=\alpha_{\{ \sigma_1,\sigma_2,\sigma_3\}}(P)$ \newline

$\bf{Example \: 3.6}$: Let $A$ be an abelian variety over $\bar{\mathbb{Q}}$, $L$ a symmetric ample line bundle on $A$. Let $f= (F_0:...:F_N): \mathbb{P}^N \rightarrow \mathbb{P}^N$ be a morphism defined by the homogeneous polynomials $F_0,..., F_N$ of same degree $d >1$ such that $0$ is the only common zero of $F_0,..., F_N.$ Set $X= A \times \mathbb{P}^N, g_1=[2] \times \mbox{id}_{\mathbb{P}^N},$ and $g_2= \mbox{id}_A \times f.$ Put $M:= p_1^*L \otimes p_2^* O_{\mathbb{P}^N}(1),$ where $p_1$ and $p_2$ are the obvious projections. Then \begin{center}
$\stackrel{(d-1) ~ \text{times}}{\overbrace{g_1^*(M) \otimes... \otimes g_1^*(M)}} \otimes g_2^*(M) \otimes g_2^*(M)  \otimes g_2^*(M) \cong M^{\otimes (4d-1)}. $
\end{center} This gives us that a canonical height $h:= \hat{h}_{\{g_1,...,g_1,g_2,g_2,g_2\}}$ exists by [10, theorem 1.2.1]. Again, if $h(P) \neq 0$, then $ \alpha_{\{g_1,...,g_1,g_2,g_2,g_2\}}(P)=\dfrac{4d-1}{d+2}$, and we can also see that $\delta_{\{g_1,...,g_1,g_2,g_2,g_2\}}= \max \{\delta_f, \delta_{[2]} \}= \max \{d, 4 \}$, which leads also to the same as the previous examples, since $\dfrac{4d-1}{d+2} <\max \{d, 4 \}.$ \newline

The next proposition is a counting orbit points result in the case of a system possibly with several maps. This result describes some information about the growth of the height counting function of the orbit of $P$ as given below.\newline

{\bf Proposition 3.7: }
\textit{Let $ P \in X_{\mathcal{F}}(\bar{K})$ whose $\mathcal{F}$-orbit is infinite, and such that the arithmetic degree $\alpha_{\mathcal{F}}(P)$ exists. Then}
\begin{center}
$\lim_{B \rightarrow \infty} \dfrac{ \# \{ n \geq 0 ; \sum_{f \in \mathcal{F}_n} h_X(f(P)) \leq B\}}{ \log B}= \dfrac{1}{ \log(k. \alpha_{\mathcal{F}}(P))} $
\end{center} \textit{and in particular,}
\begin{center}
$\lim \inf _{B \rightarrow \infty} (\# \{Q \in \mathcal{O}_{\mathcal{F}}(P); h_X(Q) \leq B \})^{\frac{1}{\log B}} \geq  k^{ \frac{1}{ \log(k. \alpha_{\mathcal{F}}(P))}}$
\end{center}

\begin{proof}
Since $\mathcal{O}_{\mathcal{F}}(P)= \infty$, it is only necessary to prove the same claim with $h_X^+$ in place of $h_X.$ For each $\epsilon >0$, there exists an $n_0(\epsilon)$ such that
\begin{center}
$(1- \epsilon) \alpha_{\mathcal{F}}(P) \leq \dfrac{1}{k} (\sum_{f \in \mathcal{F}_n} h^+_X(f(P)))^{\frac{1}{n}} \leq (1 + \epsilon)  \alpha_{\mathcal{F}}(P) $ for all $n \geq n_0(\epsilon).$
\end{center} It follows that
\begin{center}
$\{ n \geq n_0(\epsilon):  (1 + \epsilon)  \alpha_{\mathcal{F}}(P) \leq \dfrac{B^{\frac{1}{n}}}{k} \} \subset \{  n \geq n_0(\epsilon):\sum_{f \in \mathcal{F}_n} h^+_X(f(P)) \leq B \}$
\end{center} and
\begin{center}
$ \{  n \geq n_0(\epsilon):\sum_{f \in \mathcal{F}_n} h^+_X(f(P)) \leq B \} \subset \{ n \geq n_0(\epsilon):  (1 - \epsilon)  \alpha_{\mathcal{F}}(P) \leq \dfrac{B^{\frac{1}{n}}}{k} \}$
\end{center} 
Counting the number of elements in these sets yields
\begin{center}
$\dfrac{ \log B}{ \log (k (1 + \epsilon)  \alpha_{\mathcal{F}}(P))} -  n_0(\epsilon) -1 \leq \# \{ n \geq 0 : \sum_{f \in \mathcal{F}_n} h^+_X(f(P)) \leq B \}$
\end{center}
and
\begin{center}
$  \# \{ n \geq 0 : \sum_{f \in \mathcal{F}_n} h^+_X(f(P)) \leq B \} \leq \dfrac{ \log B}{ \log (k (1 - \epsilon)  \alpha_{\mathcal{F}}(P))} + n_0(\epsilon) +1 $
\end{center}
Dividing by $\log B$ and letting $B \rightarrow \infty$ gives
\begin{center}
$ \dfrac{ 1}{ \log (k (1 + \epsilon)  \alpha_{\mathcal{F}}(P))} \leq \lim \inf_{ B \rightarrow \infty}  \dfrac{\# \{ n \geq 0 : \sum_{f \in \mathcal{F}_n} h^+_X(f(P)) \leq B \} }{ \log B} $
\end{center} and
\begin{center}
$ \lim \sup_{ B \rightarrow \infty}  \dfrac{\# \{ n \geq 0 : \sum_{f \in \mathcal{F}_n} h^+_X(f(P)) \leq B \} }{ \log B} \leq  \dfrac{ 1}{ \log (k (1 - \epsilon)  \alpha_{\mathcal{F}}(P))}$
\end{center}
Since the choice for $\epsilon$ is arbitrary, and the $ \lim \inf$ is less or equal to the $\lim \sup$, this finishes the proof that 
\begin{center}
$ \lim_{ B \rightarrow \infty}  \dfrac{\# \{ n \geq 0 : \sum_{f \in \mathcal{F}_n} h^+_X(f(P)) \leq B \} }{ \log B}=  \dfrac{ 1}{ \log (k .\alpha_{\mathcal{F}}(P))}$
\end{center} Moreover, we also have that
 \begin{center}
$\{ n \geq 0 : \sum_{f \in \mathcal{F}_n} h^+_X(f(P)) \leq B \} \subset  \{ n \geq 0 :  h^+_X(f(P)) \leq B $ for all $ f \in \mathcal{F}_n \}$
\end{center} and thus
\begin{center}
$\dfrac{ \log B}{ \log (k (1 + \epsilon)  \alpha_{\mathcal{F}}(P))} -  n_0(\epsilon) -1 \leq \# \{ n \geq 0 :  h^+_X(f(P)) \leq B $ for all $ f \in \mathcal{F}_n \}$
\end{center} This implies that
\begin{center}
$\dfrac{ k^{\frac{ \log B}{ \log (k (1 + \epsilon)  \alpha_{\mathcal{F}}(P))} -  n_0(\epsilon)} -1}{k-1} \leq \# \{Q \in \mathcal{O}_{\mathcal{F}}(P); h_X^+(Q) \leq B \}$
\end{center}
Taking $\frac{1}{\log B}$-roots and letting $B \rightarrow \infty$ gives
\begin{center}
$  k^{ \frac{1}{ \log (k. \alpha_{\mathcal{F}}(P))}} \leq \lim \inf _{B \rightarrow \infty} (\# \{Q \in \mathcal{O}_{\mathcal{F}}(P); h_X^+(Q) \leq B \})^{\frac{1}{\log B}}.$
\end{center}
\end{proof}

We finish this section by stating the following elementary lemma from linear algebra. This lemma will be useful in the following sections. \newline

{\bf Lemma 3.8:} \textit{Let $A=(a_{ij}) \in M_r( \mathbb{C})$ be an $r$-by-$r$ matrix. Let $||A||=\max |a_{ij}|$, and let $\rho (A)$ denote the spectral radius of $A$. Then there are constants $c_1$ and $c_2$, depending on $A$, such that}
\begin{center}
$c_1\rho (A)^n \leq ||A^n|| \leq c_2 n^r \rho (A)^n$ for all $n \geq 0.$
\end{center}
\textit{In particular, we have $\rho (A) = \lim_{n \rightarrow \infty} ||A^n||^{ \frac{1}{n}}$. }

\begin{proof} 
See [12, lemma 14]
\end{proof}

\section{Some divisor and height inequalities for rational maps}

We let $h,g :X \dashrightarrow X$ be rational maps, and $ f \in \mathcal{F}_n$ for $\mathcal{F}=\{f_1,...,f_k\}$ a dynamical system of self-rational maps on $X$. The aim of this section is mainly to prove the next result below. It states that the action of $f \in \mathcal{F}_n$ on the vector space NS$(X)_{\mathbb{R}}$ is related with the actions of the maps $f_1,...,f_k$ by the existence of certain inequalities. This result guarantees, for instance, that the dynamical degree exists, and afterwards will also be important in order to claim and prove that $h^+_X(f(P)) \leq O(1).k^n.(\delta_{\mathcal{F}} + \epsilon)^n  h^+_X(P)$ for all $f \in \mathcal{F}_n$. 

{\bf Proposition 4.1:}
\textit{Let $X$ be a smooth projective variety, and fix a basis $D_1,..., D_r$ for the vector space NS$(X)_{\mathbb{R}}$. A dominant rational map $h: X \dashrightarrow X$ induces a linear map on  NS$(X)_{\mathbb{R}}$, and we write}
\begin{center}
$h^*D_j \equiv \sum_{i=1}^r a_{ij}(h)D_i$ \textit{and} $A(h)=(a_{ij}(h)) \in M_r(\mathbb{R}).$
\end{center}
\textit{We let $||.||$ denote the sup norm on $M_r(\mathbb{R}).$ Then there is a constant $C \geq 1$ depending on $D_1,...,D_r$ such that for any dominant rational maps $h,g :X \dashrightarrow X,$ any $n \geq 1$, and any $ f \in \mathcal{F}_n$ we have}
\begin{center}
$\quad \quad \quad|| A( g \circ h)|| \leq C ||A(g)|| . ||A(h)||$ \newline 
$||A(f)|| \leq C.(r . \max_{i=1,...,k.}||A(f_i)||)^n.$
\end{center}

The proof of this result will be made in the sequel. An immediate corollary of this is the convergence of the limit defining the dynamical degree. \newline

{\bf Corollary 4.2:}
\textit{The limit $\delta_{\mathcal{F}}=\lim \sup_{n\rightarrow \infty}\rho(\mathcal{F}_n)^{\frac{1}{n}}$ exists.}

\begin{proof}
With notation as in the statement of proposition 4.1, we have
\begin{center}
$\rho(\mathcal{F}_n)= \max_{f \in \mathcal{F}_n} \rho (f^*,$NS$(X)_{\mathbb{R}})= \max_{f \in \mathcal{F}_n} \rho (A(f))$
\end{center} 

Denoting $||A(\mathcal{G})|| = \max_{g \in \mathcal{G}}||A(g)||$, where $\rho(\mathcal{G}):= \rho(g)$ for $\mathcal{G}$ dynamical system and $g \in \mathcal{G}$, proposition 4.1 give us that
\begin{center}
$\log ||A(\mathcal{F}_{n+m})||  \leq \log ||A(\mathcal{F}_{m})|| + \log ||A(\mathcal{F}_{n})||+ O(1)$
\end{center} Using this convexity estimate, we can see that $\dfrac{1}{n}\log ||A(\mathcal{F}_n)||$ converges. Indeed, if a sequence $(d_n)_{n \in \mathbb{N}}$ of nonnnegative real numbers satisfies $d_{i+j} \leq d_i + d_j$, then after fixing a integer $m$ and writing $n=mq+r$ with $0 \leq r \leq m-1,$ we have
\begin{center}
$\dfrac{d_n}{n} = \dfrac{d_{mq+r}}{n} \leq  \dfrac{(qd_m+d_r)}{n}= \dfrac{d_m}{m} \dfrac{1}{(1+ r/mq)} + \dfrac{d_r}{n} \leq \dfrac{d_m}{m} + \dfrac{d_r}{n}.$
\end{center} Now take the limsup as $n \rightarrow \infty$, keeping in mind that $m$ is fixed and \newline $ r \leq m-1$, so $d_r$ is bounded. This gives
\begin{center}
$\lim \sup_{n \rightarrow \infty} \dfrac{d_n}{n} \leq \dfrac{d_m}{m}.$
\end{center} taking the infimum over $m$ shows that
\begin{center}
$\lim \sup_{ n \rightarrow \infty}  \dfrac{d_n}{n} \leq \inf_{ m \geq 1}  \dfrac{d_m}{m} \leq \lim \inf_{ m \rightarrow \infty}  \dfrac{d_m}{m} ,$
\end{center} and hence all three quantities must be equal.

As the sequence $(||A(\mathcal{F}_{n})||^{1/n})_{n \in \mathbb{N}}$ is convergent and therefore bounded, lemma 3.8 guarantees that the sequence $(\rho(\mathcal{F}_n)^{1/n})_{n \in \mathbb{N}}$ is bounded as well.
\end{proof}

We also conjecture that the limit $\lim_{n\rightarrow \infty}\rho(\mathcal{F}_n)^{\frac{1}{n}}$ exists and is a birational invariant. The proof for dynamical degrees of systems with only one map given in [6, prop. 1.2] should be extented naturally for our present definition of degree with several maps. In the mentioned article, the dynamical degree is firstly defined using currents, and afterwards such definition is proved to coincide with the one using the limit of roots of spectral radius. Such result can be worked out in some future paper. Thus, from now on, we assume that \begin{center} $\delta_{\mathcal{F}}:=\lim_{n\rightarrow \infty}\rho(\mathcal{F}_n)^{\frac{1}{n}}$, \end{center} and that it exists.

We start the proof of proposition 4.1 stating the following auxiliar proposition and lemmas whose proofs can be found in [12]:\newline

{\bf Proposition 4.3:}
\textit{Let $X^{(0)}, X^{(1)}, ..., X^{(m)}$ be smooth projective varieties of the same dimension $N$, and let $f^{(i)}: X^{(i)} \dashrightarrow X^{(i-1)}$ be dominant rational maps for $1 \leq i \leq m.$ Let $D$ be a nef divisor on $X^{(0)}$. Then for any nef divisor $H$ on $X^{(m)}$, we have}
\begin{center}
$(f^{(1)} \circ f^{(2)} \circ ... \circ f^{(m)})^*D . H ^{N-1} \leq (f^{(m)})^*... (f^{(2)})^* (f^{(1)})^*D.H^{N-1}.$
\end{center}

\begin{proof} See [12, Prop. 17] \end{proof}
For the lemmas, we need to set the following notation: \newline \begin{itemize}
\item $N:$  The dimension of $X$, wich we assume is at least 2. \newline
\item $\mbox{Amp}(X)$: The ample cone in NS$(X)_{\mathbb{R}}$ of all ample $\mathbb{R}-$divisors. \newline
\item $\mbox{Nef}(X)$:  The nef cone in NS$(X)_{\mathbb{R}}$ of all nef $\mathbb{R}-$divisors. \newline
\item $\mbox{Eff}(X)$:  The effective cone in NS$(X)_{\mathbb{R}}$ of all effective $\mathbb{R}-$divisors. \newline
\item $\overline{\mbox{Eff}}(X): $ The $\mathbb{R}-$closure of Eff$(X)$. \newline
\end{itemize}

As described in [5, section 1.4], we have the facts
\begin{center}
Nef$(X)=\overline{\mbox{Amp}}(X)$ and Amp$(X)=$ int$($Nef$(X)).$
\end{center} In particular,  since Amp$(X) \subset $ Eff$(X)$, it follows that  Nef$(X) \subset \overline{\mbox{Eff}}(X).$\newline

{\bf Lemma 4.4:} \textit{With notation as above, let $D \in  \overline{\mbox{Eff}}(X) - \{0\}$ and $H \in $ Amp$(X).$ Then }
\begin{center}
$D.H^{N-1} > 0.$ 
\end{center}

\begin{proof}
See [12, lemma 18]
\end{proof} \newpage
{\bf Lemma 4.5:}
\textit{Let $H \in $ Amp$(X)$, and fix some norm $|.|$ on the $\mathbb{R}-$vector space NS$(X)_{\mathbb{R}}$. Then there are constants $C_1, C_2 >  0$ such that}
\begin{center}
$C_1|v| \leq v.H^{N-1} \leq C_2|v|$ for all $v \in  \overline{\mbox{Eff}}(X).$
\end{center}
\begin{proof}
See [12, lemma 19]
\end{proof}
Now we start the proof of proposition 4.1. We fix a norm $|.|$ on the $\mathbb{R}-$vector space NS$(X)_{\mathbb{R}}$ as before. Additionally, for any $A:$ NS$(X)_{\mathbb{R}} \rightarrow $ NS$(X)_{\mathbb{R}}$ linear transformation, we set
\begin{center}
$||A||^{\prime} = \sup_{ v \in \mbox{Nef} - \{0\}} \dfrac{|Av|}{|v|},$
\end{center} which exists because the set $\overline{\mbox{Eff}}(X) \cap \{ w \in \mbox{NS}(X)_{\mathbb{R}} : |w| = 1 \}$ is compact. \newline
We note that for linear maps $A,B \in $ End(NS$(X)_{\mathbb{R}})$ and $c \in \mathbb{R}$ we have
\begin{center}
$|| A + B||^{\prime} \leq ||A||^{\prime} + ||B||^{\prime}$ and $||cA||^{\prime}=|c|||A||^{\prime}.$
\end{center}
Further, since Nef$(X)$ generates NS$(X)_{\mathbb{R}}$ as an $\mathbb{R}-$vector space, we have $||A||^{\prime}=0$ if and only if $A=0.$ Thus $||.||^{\prime}$ is an $\mathbb{R}-$norm on End(NS$(X)_{\mathbb{R}}).$ \newline
Similarly, for any linear map $A:$ NS$(X)_{\mathbb{R}} \rightarrow $NS$(X)_{\mathbb{R}},$ we set
\begin{center}
$||A||^{ \prime \prime}= \sup_{ v \in \mbox{Eff} - \{0\}} \dfrac{|Aw|}{|w|},$
\end{center} then $||.||^{ \prime \prime}$ is an $\mathbb{R}-$norm on  End(NS$(X)_{\mathbb{R}}).$

We note that $ \overline{\mbox{Eff}}(X)$ is preserved by $f^*$ for $f$ self-rational map on $X$, and that Nef$(X) \subset  \overline{\mbox{Eff}}(X).$ Thus if $v \in $Nef$(X),$ then $ g^*v$ and $h^*v$ belong to  $ \overline{\mbox{Eff}}(X).$ This allows us to compute \newline \newline
$||(g \circ h)^*||^{\prime}=\sup_{ v \in \mbox{Nef}(X) - \{0\}} \dfrac{|(g \circ h)^*v|}{|v|}~~~~~~~~~~~~~~~~~~~~~~ ~~~~~~~~~~~ ~~~~~~~~~~~~~\newline \newline
\leq C_1^{-1}\sup_{ v \in \mbox{Nef}(X) - \{0\}} \dfrac{(g \circ h)^*v.H^{N-1}}{|v|} $ from lemma 4.5 ~~~~~~~~~~~~~~~~~~~~~~~~ \newline \newline 
 $ \leq  C_1^{-1}\sup_{ v \in \mbox{Nef}(X) - \{0\}} \dfrac{(h^*  g^*v).H^{N-1}}{|v|} $ from proposition 4.3 ~~~~~~~~~~~~~~~~~~~~\newline \newline
$ =  C_1^{-1}\sup_{ v \in \mbox{Nef}(X) - \{0\}, g^*v \neq 0} \dfrac{(h^*  g^*v).H^{N-1}}{|v|} ~~~~~~~~~~~~~~~~~~~~~~~~~~~~~~~~~~~~~~~~~~~\newline \newline
=  C_1^{-1}(\sup_{ v \in \mbox{Nef}(X) - \{0\}, g^*v \neq 0} \dfrac{(h^*  g^*v).H^{N-1}}{|g^*v|} . \dfrac{|g^*v|}{|v|})~~~~~~~~~~~~~~~~~~~~~~~~~~~~~~~~~~ \newline \newline
 \leq  C_1^{-1}(\sup_{ v \in \mbox{Nef}(X) - \{0\}, g^*v \neq 0} \dfrac{(h^*  g^*v).H^{N-1}}{|g^*v|}) .(\sup_{ v \in \mbox{Nef} - \{0\}} \dfrac{|g^*v|}{|v|})~~~~~~~~~~~~~~ \newline \newline%
=  C_1^{-1}(\sup_{ v \in \mbox{Nef}(X) - \{0\}, g^*v \neq 0} \dfrac{(h^*  g^*v).H^{N-1}}{|g^*v|}) . || g^*||^{\prime}~~~~~~~~~~~~~~~~~~~~~~~~~~~~~~~~~~ \newline \newline
\leq  C_1^{-1}(\sup_{ w \in \overline{\mbox{Eff}}(X) - \{0\}} \dfrac{ (h^*w).H^{N-1}}{|w|}) . || g^*||^{\prime}$  since $ g^*v \in  \overline{\mbox{Eff}}(X)~~~ ~~~~~~~~~~~~~~\newline \newline
\leq C_1^{-1} C_2 (\sup_{ w \in \overline{\mbox{Eff}}(X) - \{0\}} \dfrac{ |h^*w|}{|w|}) . || g^*||^{\prime} $ from lemma 4.5~~~~~~~~~~~~~~~~~~~~~~~~~ \newline \newline 
$=  C_1^{-1} C_2 ||h^*||^{\prime \prime}. ||g^*||^{\prime}.~~~~~~~~~~~~~~~~~~~~~~~~~~~~~~~~~~~~~~~~~~~~~~~~~~~~~~~~~~~~~~~~~~~~~~~~~~~~$ \newline

We remember that we defined $||.||$ to be the sup norm on $M_r(\mathbb{R})=$ End$($NS$(X)_{\mathbb{R}}$, where the identification is via the given basis $D_1,..., D_r$ of  NS$(X)_{\mathbb{R}}$. We thus have three norms $||.||, ||.||^{\prime}$ and $||.||^{\prime \prime}$ on End$($NS$(X)_{\mathbb{R}}$, so there are positive constants $C_3^{\prime}, C_4^{\prime}, C_3^{\prime \prime}$ and $ C_4^{\prime \prime}$ such that 
\begin{center}
$C_3^{\prime}|| \gamma|| \leq || \gamma||^{\prime} \leq C_4^{\prime}|| \gamma||$ and $C_3^{\prime \prime}|| \gamma|| \leq || \gamma||^{\prime \prime} \leq C_4^{\prime \prime}|| \gamma||$ l $\forall \gamma \in $  End$($NS$(X)_{\mathbb{R}}.$
\end{center}
 Hence \newline \newline
$||A(g \circ h)||=||(g \circ h)^*|| \leq C_3^{\prime -1}||(g \circ h)^*||^{\prime}~~~~~~~ \newline \newline
\leq C_3^{\prime -1} C_1^{-1} C_2  ||h^*||^{\prime \prime}. ||g^*||^{\prime} ~~~~~~~~~~~~~~~~~~~~~~~~~~~~~\newline \newline
\leq  C_3^{\prime -1} C_1^{-1} C_2  C_4^{\prime}  C_4^{\prime \prime}  ||h^*||. ||g^*|| ~~~~~~~~~~~~~~~~~~~~~~~~\newline \newline
= C_3^{\prime -1} C_1^{-1} C_2  C_4^{\prime}  C_4^{\prime \prime}  ||A(h)||. ||A(g)||.~~~~~~~~~~~~~~~~~~~~~~~~~~ ~~~~~~~~$ 
\newline \newline
Similarly, if $v \in $ Nef$(X), f := f_{i_1} \circ ... \circ f_{i_n}\in \mathcal{F}_n$, then $f^*v \in \overline{\mbox{Eff}}(X).$ A similar calculation gives \newline \newline
$||f^*||^{\prime}=\sup_{ v \in \mbox{Nef}(X) - \{0\}} \dfrac{|f^*v|}{|v|} ~~~~~~~~~~~~~~~~~~~~~~~~~     ~~~~~~~~~~~~~~~~~~~~~~~~~~~~~~~~~~~~~~~~~~~~\newline \newline
\leq C_1^{-1}\sup_{ v \in \mbox{Nef}(X) - \{0\}} \dfrac{(f^*v).H^{N-1}}{|v|} $ from lemma 4.5 ~~~~~~~~~~~~~~~~~~~~~~~~~~~~~~~~~~~~~~\newline \newline 
 $=  C_1^{-1}\sup_{ v \in \mbox{Nef}(X) - \{0\}} \dfrac{( f_{i_1} \circ ... \circ f_{i_n})^*v.H^{N-1}}{|v|} ~~~~~~~~~~~~~~~~~~~~~~~~~~~~~~ ~~~~~~~~~~~~~~~~~\newline \newline
\leq  C_1^{-1}\sup_{ v \in \mbox{Nef}(X) - \{0\}} \dfrac{((f_{i_n})^*...  (f_{i_1})^*v).H^{N-1}}{|v|} $ from proposition 4.3~~~~~~~~~~~~~~~~~\newline \newline 
$ \leq   C_1^{-1} C_2 (\sup_{ v \in \mbox{Nef}(X) - \{0\}} \dfrac{ |(f_{i_n})^*...  (f_{i_1})^*v|}{|v|})  $ from lemma 4.5 ~~~~~~~~~~~ ~~~~~~~~~~~~~~\newline \newline 
$=  C_1^{-1} C_2. ||(f_{i_n})^*...  (f_{i_1})^*||^{\prime}.~~~~~~~~~~~~~~~~~~~~~~~~~~~~~~~~~~~~~~~~~~~~~~~ ~~~ ~~~~~~~~~~~~~~~~~~~~~~~~$ \newline\newline  Hence \newline \newline
$||A(f)||=||f^*|| \leq  C_3^{\prime -1}||f^*||^{\prime}~~~~~~~~~~~~~~~~~~~~~~  ~~~~~~~~~~ \newline \newline
\leq C_3^{\prime -1} C_1^{-1} C_2  ||(f_{i_n})^*...  (f_{i_1})^*||^{\prime}  ~~~~~~~~~~~~~~~~~~~~~  ~~~~~~~~~~\newline \newline
\leq  C_3^{\prime -1} C_1^{-1} C_2  C_4^{\prime}  C_4^{\prime \prime} ||(f_{i_n})^*...  (f_{i_1})^*|| ~~~~~~~~~~~~~   ~~~~~~~~~~~\newline \newline
\leq  C_3^{\prime -1} C_1^{-1} C_2  C_4^{\prime}  C_4^{\prime \prime} r^n ||(f_{i_n})^*||...  ||(f_{i_1})^*||~~~~~~~  ~~~~~~~~~~~ \newline \newline
\leq C_3^{\prime -1} C_1^{-1} C_2  C_4^{\prime}  C_4^{\prime \prime}.[r. \max_{i=1,...,k.}||A(f_i)||]^n$, ~~~~~~~~~~~~~~~~~~~~~~~~~~~\newline \newline
As we wanted to show. \newline

As it was said in the beginning of this section, the next proposition is a height inequality for rational maps, with eyes towards future applications. \newline

{\bf Proposition 4.6:}
\textit{Let $X/\bar{K}$ and $ Y/\bar{K}$ be smooth projective varieties, \newline let $f:Y \dashrightarrow X$ be a dominant rational map defined over $\bar{K}$, let $D \in$ Div$(X)$ be an ample divisor, and fix Weil height functions $h_{X,D}$ and $h_{Y,f^*D}(P)$ associated to $D$ and $f^*D.$ Then }
\begin{center}
$h_{X,D} \circ f(P) \leq h_{Y,f^*D}(P) + O(1)$ \textit{for all} $P \in (Y - I_f)(\bar{K}),$
\end{center} \textit{where the $O(1)$ bound depends on $X,Y, f,$ and the choice of height functions, but is independent of $P$.}

\begin{proof}
See [12, Prop. 21].
\end{proof}


\section{A bound for the sum of heights on iterates} 
This section is devoted for the proof of a quantitative upper bound for $\sum_{f \in \mathcal{F}_n} h^+_X(f(P))$ in terms of the dynamical degree $\delta_{\mathcal{F}}$ of the system.This is one of the main results of this work, and is stated below. As a corollary, we see that the arithmetic degree of any point is upper bounded by the dynamical degree of the system. \newline

{\bf Theorem 5.1:}
\textit{Let $K$ be a number field or a one variable function field of characteristic $0$ , let $\mathcal{F}=\{f_1,...,f_k\}$ be a set of dominant self rational maps on $X$  defined over $K$ as stated before, let $h_X$ be a Weil height on $X(\bar{K})$ relative to an ample divisor, let $h^+_X= \max \{h_X, 1 \}$, and let $\epsilon >0$. Then there exists a positive constant $C=C(X, h_X, f, \epsilon)$ such that for all $P \in X_{\mathcal{F}}(\bar{K})$ and all $n \geq 0$,}
\begin{center}
$\sum_{f \in \mathcal{F}_n} h^+_X(f(P)) \leq C. k^n.(\delta_{\mathcal{F}} + \epsilon)^n .  h^+_X(P).$
\end{center} \textit{In particular, $h^+_X(f(P)) \leq  C. k^n.(\delta_{\mathcal{F}} + \epsilon)^n .  h^+_X(P)$ for all $f \in \mathcal{F}_n.$} \newline \newline
 Before proving the theorem, we note that it implies the fundamental inequality $\bar{\alpha}_{\mathcal{F}}(P) \leq \delta_{\mathcal{F}}.$\newline

{\bf Corollary 5.2:}
\textit{Let $P \in  X_{\mathcal{F}}(\bar{K}).$ Then}
\begin{center}
$\bar{\alpha}_{\mathcal{F}}(P) \leq \delta_{\mathcal{F}}.$
\end{center} 

\begin{proof} Let $ \epsilon >0.$ Then
\begin{center}
$\quad \quad\bar{\alpha}_{\mathcal{F}}(P) = \frac{1}{k} \lim \sup_{n \rightarrow \infty} \{\sum_{f \in \mathcal{F}_n} h^+_X(f(P))\}^{\frac{1}{n}}  $ by definition of $\bar{\alpha}_{\mathcal{F}}~ ~~~~~~~\newline \newline
\leq  \lim \sup_{n \rightarrow \infty} ( C. (\delta_{\mathcal{F}} + \epsilon)^n .  h^+_X(P))^{\frac{1}{n}}  $ from theorem 5.1 ~~~~~~~~~~~~~~~ ~~~\newline \newline
$ =\delta_{\mathcal{F}} + \epsilon.\quad \quad \quad \quad \quad \quad\quad \quad\quad \quad\quad \quad\quad \quad\quad \quad\quad \quad\quad \quad\quad \quad\quad \quad  ~~~~~~~~~~~~~~~~~~~~~~~~~~~~~~~~~~~~~~~~~~~~~~~~~~~~~~~~~~~~~~ ~~~~~~~~~~~~~~~~~~~~~~~~~$
\end{center}
This holds for all $\epsilon>0$, which proves that  $\bar{\alpha}_{\mathcal{F}}(P) \leq \delta_{\mathcal{F}}.$
\end{proof} 

{\bf Lemma 5.3:}
\textit{Let $E \in$ Div$(X)_{\mathbb{R}}$ be a divisor that is algebraic equivalent to 0, and fix a height function $h_E$ associated to $E.$ Then there is a constant $C=C(h_X, h_E)$ such that}
\begin{center}
$|h_E(P)| \leq C \sqrt{h_X^+(P)} $ \textit{for all} $P \in  X(\bar{K}).$
\end{center}

\begin{proof}
See for example the book of Diophantine Geometry of Hindry-Silverman[8, Theorem B.5.9].
\end{proof}
Theorem 5.1 will be a consequence from the slightly weaker result:\newline

{\bf Theorem 5.4:}
\textit{Let $K$ be a number field or a one variable function field of characteristic $0$ , let $\mathcal{F}=\{f_1,...,f_k\}$ be a set of dominant self rational maps on $X$  defined over $K$, let $h_X$ be a Weil height on $X(\bar{K})$ relative to an ample divisor, let $h^+_X= \max \{h_X, 1 \}$, and let $\epsilon >0$. Then there exists a positive constant $C=C(X, h_X, f, \epsilon)$, and $t$ positive integer such that for all $P \in X_{\mathcal{F}}(\bar{K})$ and all $n \geq 0$,}
\begin{center}
$\sum_{f \in \mathcal{F}_{nt}} h^+_X(f(P)) \leq C. k^{nt}.(\delta_{\mathcal{F}} + \epsilon)^{nt} .  h^+_X(P).$
\end{center} 

Before proving it and then deduce theorem 5.1, we state and prove two auxiliar short lemmas.\newline

{\bf Lemma 5.5:} \textit{In the situation above, there is a constant $C \geq 1$ such that
\begin{center}
$\sum_{f \in \mathcal{F}_n} h^+_X(f(P)) \leq k^n.C^n .  h^+_X(P).$
\end{center} 
for all $P \in X_{\mathcal{F}}(\bar{K})$.}
\begin{proof}
We take $H$ an ample divisor on $X$, $h_H \geq 1$ and $ h_{f_i^*H}$ height functions associated to $H$ and $f_i^*H$ respectively, so that \begin{center}  $h_H(f_i(P)) \leq h_{f_i^*H}(P) + O(1) $ \end{center}
 for all $P \in X_{\mathcal{F}}(\bar{K})$ , with $O(1)$ depending on $H, f_i, f_i^*H, h_H, h_{f_i^*H},$ but not on $P$. Then, for $C$ enough large,  we find that $h_{f_i^*H}(P) + O(1) \leq C h_H(P)$, and so $h_H(f_i(P)) \leq C h_H(P)$ for all $P \in X_{\mathcal{F}}(\bar{K})$, which yields \begin{center}$\sum_{f \in \mathcal{F}_n} h_H(f(P)) \leq k^n.C^n .  h_H(P).$ \end{center}
 The proof is finished since $h_H$ and $h_X$ are associated with ample divisors, and therefore are commensurate.
\end{proof}

{\bf Lemma 5.6:} \textit{Let $\mathcal{A}_0:=\{a_0\}, a_0 \geq 1$, $k$ fixed, and for each $l \in \mathbb{N}$, $\mathcal{A}_l$ a set with $k^l$ positive  real numbers such that
\begin{center} $\sum_{a \in \mathcal{A}_n} a \leq \sum_{a \in \mathcal{A}_{n-1}} a + C_1(\sum_{a \in \mathcal{A}_{n-1}} \sqrt{a} +  \sum_{a \in \mathcal{A}_{n-1}} \sqrt{a + C_2})$ for all $n \geq 1,$ \end{center} where $C_1, C_2$ are non-negative constants. Then there exists a positive constant $C$ depending only on $C_1, C_2$ such that
\begin{center}
$\sum_{a \in \mathcal{A}_n} a \leq k^{n-1}.C.n^2. a_0$
\end{center}}
\begin{proof}
$\sum_{a \in \mathcal{A}_n} a \leq \sum_{a \in \mathcal{A}_{n-1}} a + C_1(\sum_{a \in \mathcal{A}_{n-1}} \sqrt{a} +  \sum_{a \in \mathcal{A}_{n-1}} \sqrt{a + C_2}) \newline 
=\sum_{a \in \mathcal{A}_{n-1}} [a + C_1( \sqrt{a} +   \sqrt{a + C_2})] \leq \sum_{a \in \mathcal{A}_{n-1}} [a + C_1\sqrt{a}( 1 +   \sqrt{1 + \dfrac{C_2}{a}})] \newline 
\leq  \sum_{a \in \mathcal{A}_{n-1}} [a + C_1\sqrt{a}( 1 +   \sqrt{1 + C_2})]= \sum_{a \in \mathcal{A}_{n-1}} [a + C_3\sqrt{a} ]$  with \newline $C_3:= C_1(1+ \sqrt{1+ C_2})$.\newline 
Thus we have $\sum_{a \in \mathcal{A}_1}a \leq \sum_{a \in \mathcal{A}_{0}} [a + C_3\sqrt{a} ]=a_0 + C_3\sqrt{a_0} \leq a_0(1 + C_3) \newline \leq a_0.C = a_0.C.k^0,$ where $C:=\max \{ \dfrac{C_3.k}{4}, 1 + C_3\}$, and we want to prove by induction that $\sum_{a \in \mathcal{A}_n}a \leq Ck^{n-1}n^2a_0.$ So we compute \newline \newline
$\sum_{a \in \mathcal{A}_{n+1}}a \leq  \sum_{a \in \mathcal{A}_{n}} [a + C_3\sqrt{a}] \leq \sum_{a \in \mathcal{A}_{n}} a+ C_3\sum_{a \in \mathcal{A}_{n}} \sqrt{a} \newline \leq \sum_{a \in \mathcal{A}_{n}}a + k^{n/2}C_3^{1/2}\sqrt{\sum_{a \in \mathcal{A}_{n}}a} \leq k^{n-1}Cn^2a_0 + k^{n/2}C_3^{1/2}\sqrt{k^{n-1}Cn^2a_0} \newline \leq k^{n-1}Cn^2a_0+ 2k^{(n-1)/2}\sqrt{C}\sqrt{k^{n-1}Cn^2a_0} \leq k^{n-1}Ca_0(n^2+\dfrac{2n}{\sqrt{a_0}})\newline \leq k^{n-1}Ca_0(n^2+2n)  \leq k^{n}\tilde{C}a_0(n+1)^2,$ where we can make $\tilde{C}:=\max \{C_3/4, 1+ C_3 \}$
\end{proof}
Now we start the proof of theorem 5.4 \newline

\textit{ Proof of theorem 5.4:} We take $D_1,..,D_r$ very ample divisors forming a basis of NS$(X)_{\mathbb{R}}$, and $H \equiv \sum c_i D_i$ ample with $c_i \geq 0$ such that $H+ D_i, H-D_i$ are all ample.

We consider a resolution of indeterminacy $p:Y \rightarrow X$ as a sequence of blowing ups working for each $f_i$, such that $g_i:=f_i \circ p$ is a morphism for each $i \leq k$, and Exc$(p)$ is the exceptional locus of $p$. For each $j\leq k, i \leq r$, we take effective divisors $\tilde{D}_i^{(j)}$ on $X$ with $\tilde{D}_i^{(j)}$ linearly equivalent to $D_i$, and such that none of the components of $g_j^*\tilde{D}_i^{(j)}$ are containded in Exc$(p)$. The divisor $Z_i^{(j)}:=p^*p_*g_j^*\tilde{D}_i^{(j)}- g_j^*\tilde{D}_i^{(j)}$ on $Y$ is effective and has support contained in Exc$(p)$. We denote $F_i^{(j)}:=g_j^*D_i$ for $i=1,...,r$, and take divisors $F_{r+1}^{(j)},...,F_s^{(j)}$ so that $F_1^{(j)},..., F_s^{(j)}$ form a basis for NS$(Y)_{\mathbb{R}}$. For $i \leq r$, we can see that $p^*p_*F_i^{(j)}-F_i^{(j)}$ and $Z_i^{(j)}$ are linearly equivalent. By [7, prop. 7.10], we can find $H^{\prime} \in$ Div$(Y)_{K}$ ample so that $p^*H - H^{\prime}$ is effective with support contained in Exc$(p)$.

We consider $g_j^*D_i \equiv \sum_{m \leq s } a^{(j)}_{mi}F_m^{(j)}$ for $i=1,...,r$ and $A^{(j)}:=(a^{(j)}_{mi})_{m,i}$ the correspondent $s \times r-$matrix. We also denote $p_*F_i^{(j)} \equiv \sum_{l \leq r } b^{(j)}_{li}D_l; i=1,...,s,$ and  $B^{(j)}:=(b^{(j)}_{li})_{l,i}$ the correspondent $r \times s-$matrix.

We see that $B^{(j)}A^{(j)}$ is a matrix representing $f_j^*$ with respect to the basis $D_1,...,D_r$. 

Let us fix some notation:

$\vec{D}:=(D_1,...,D_r), \vec{F}^{(j)}:=(F_1^{(j)},..., F_s^{(j)}), \vec{Z}^{(j)}:=(Z_1^{(j)},..., Z_s^{(j)}), \vec{c}:=(c_1,...,c_r), \newline E^{(j)}:=g_j^*H-<A^{(j)} \vec{c}, \vec{F}^{(j)}>, \vec{E^{\prime}}^{(j)}= ({E^{\prime}_1}^{(j)},...,{E^{\prime}_s}^{(j)}):=p_*\vec{F}^{(j)} - {B^{(j)}}^T\vec{D}. $ We note that $E^{(j)}$ and $\vec{E^{\prime}}^{(j)}$ are numerically zero divisors for each $j$.

We choose height functions $h_{D_1},..., h_{D_r}$ for $D_1,...D_r$ respectively, and $h_H \geq 1$ with respect to $H$ such that $h_H \geq |h_{D_i}|$ for each $i \leq r$. All of these functions are independent of $\mathcal{F}$. Defining $h_{F^{(j)}_i}:=h_{D_i}\circ g_j, i=1,...,r$ height functions associated with $F_i^{(j)}$. For $i=r+1,...,s$ we fix height functions $h_{p_*F^{(j)}_i}$ with respect to the divisors $p_*F^{(j)}_i$, and we denote: $h_{\vec{D}}:=(D_1,...,D_r), h_{\vec{F}^{(j)}}:=(h_{F_1^{(j)}},...,h_{F_s^{(j)}}), \newline h_{p_*\vec{F}^{(j)}}:=(h_{p_*F_1^{(j)}},...,h_{p_*F_s^{(j)}}), h_{\vec{E^{\prime}}^{(j)}}:=(h_{{E^{\prime}_1}^{(j)}},...,h_{{E^{\prime}_s}^{(j)}}) = h_{p_*\vec{F}^{(j)}}-{B^{(j)}}^Th_{\vec{D}}, \newline h_{\vec{Z}^{(j)}}:= (h_{Z_1^{(j)}},...,h_{Z_s^{(j)}})=h_{p_*\vec{F}^{(j)}}\circ p - h_{\vec{F}^{(j)}},$ where  $h_{\vec{Z}_i^{(j)}}$ and $h_{\vec{E^{\prime}}_i^{(j)}}$ are height functions associated with the divisors $\vec{Z}_i^{(j)}$ and $\vec{E^{\prime}}_i^{(j)}$. Also, define \begin{center} $h_{E^{(j)}}:=h_H \circ g_j - <A^{(j)} \vec{c}, h_{\vec{F^{(j)}}}>$. \end{center} We can suppose that $h_{Z_i^{(j)}} \geq 0$ on $Y-Z_i^{(j)}$. We can fix a height function $h_{H^{\prime}} \geq 1$ related to $H^{\prime}$, and a height function $h_{p^*H - H^{\prime}}$ related to $p^*H - H^{\prime}$ satisfying $h_{p^*H - H^{\prime}} \geq 0$ on $Y-$Exc$(p)$. Since $E^{(j)}$ and ${E_i^{\prime}}^{(j)}$ are numerically equivalent to zero, there exists a positive constant such that $|h_{{E}^{(j)}}| \leq C \sqrt{h_{H^{\prime}}}$ and $|h_{{E_i^{\prime}}^{(j)}}|\leq C \sqrt{h_H}$. Also, there exists a constant $\gamma \geq 0$ such that $h_H \circ p \geq h_{p^*H - H^{\prime}} + h_{H^{\prime}}- \gamma$ on $Y(\bar{K})$. Finally, if we denote by $M(f_j)$ the matrix representing $f_j^*$, linear map on NS$(X)_{\mathbb{R}},$ with respect to the basis $D_1,...,D_r$, $||M(f_j)||$ the maximum absolute value of its coefficients (norm of a matrix), then we make the notation $||M(\mathcal{F}_n)||:= \max_{f \in \mathcal{F}_n}||M(f)||$. 

For $P \in X_{\mathcal{F}}(\bar{K}), n \geq 1$, we compute: \newline \newline
$\sum_{f \in \mathcal{F}_n} h_H(f(P))= \sum_{i \leq k}\sum_{f \in \mathcal{F}_{n-1}}h_H(f_i(f(P))) \newline \newline = \sum_{i \leq k,f \in \mathcal{F}_{n-1}}[(h_H \circ g_i)(p^{-1}f(P))-<A^{(i)}\vec{c}, h_{p_*F^{(i)}}\circ p>(p^{-1}f(P))\newline \newline +<A^{(i)}\vec{c}, h_{p_*F^{(i)}}>(f(P))]=\sum_{i\leq k,f \in \mathcal{F}_{n-1}}[<A^{(i)}\vec{c}, h_{F^{(i)}}-h_{p_*F^{(i)}}\circ p>(p^{-1}f(P))\newline \newline + h_{E^{(i)}}(p^{-1}f(P))+ <B^{(i)}A^{(i)}\vec{c},h_{\vec{D}}>(f(P)) +<A^{(i)}\vec{c},h_{{E^{\prime}}^{(i)}}>(f(P))]\newline \newline =\sum_{i\leq k,f \in \mathcal{F}_{n-1}}[<\vec{c},-h_{Z^{(i)}}>(p^{-1}f(P))+h_{E^{(i)}}(p^{-1}f(P))\newline \newline + <B^{(i)}A^{(i)}\vec{c},h_{\vec{D}}>(f(P))+ <\vec{c},{A^{(i)}}^Th_{{E^{\prime}}^{(i)}}>(f(P))] \newline \newline \leq \sum_{i\leq k,f \in \mathcal{F}_{n-1}}[h_{E^{(i)}}(p^{-1}f(P)) +<B^{(i)}A^{(i)}\vec{c},h_{\vec{D}}>(f(P))\newline \newline +<\vec{c},{A^{(i)}}^Th_{{E^{\prime}}^{(i)}}>(f(P))]\leq  \sum_{i\leq k,f \in \mathcal{F}_{n-1}}[r^2||\vec{c} ||||B^{(i)}A^{(i)}|| h_H(f(P))\newline \newline+r||\vec{c}||C\sqrt{h_H(f(P))} +C \sqrt{h_{H^{\prime}}(p^{-1}f(P))}]\newline \newline \leq \sum_{i\leq k,f \in \mathcal{F}_{n-1}}[r^2||\vec{c} ||||B^{(i)}A^{(i)}|| h_H(f(P))+r||\vec{c}||C\sqrt{h_H(f(P))} +C \sqrt{h_{H^{\prime}}(f(P)) + \gamma}]$,\newline \newline where the last follows because $h_H \circ p \geq h_{p^*H - H^{\prime}} + h_{H^{\prime}}- \gamma$ on $Y(\bar{K})$ and $h_{p^*H - H^{\prime}} \geq 0$ on $Y-$Exc$(p)$.\newline

Denoting by $R:=\max_i \{1, r^2||\vec{c} ||||B^{(i)}||||A^{(i)}||\}$, and dividing the whole inequality above by $R^n$, we obtain \newline \newline
$\dfrac{1}{R^{n}}\sum_{f \in \mathcal{F}_n} h_H(f(P))\newline \leq k.[\sum_{f \in \mathcal{F}_{n-1}} \dfrac{h_H(f(P))}{R^n}+r||\vec{c}||C \sum_{f \in \mathcal{F}_{n-1}} \sqrt{\dfrac{h_H(f(P))}{R^{n-1}}}+C\sum_{f \in \mathcal{F}_{n-1}} \sqrt{\dfrac{h_H(f(P))}{R^{n-1}}+ \gamma}],$ \newline \newline
 which, by lemma 5.6, implies that \begin{center}
$\sum_{f \in \mathcal{F}_n} h_H(f(P))\leq C_1 k^n n^2 R^n h_H(P),$
\end{center} for a positive constant $C_1.$

Fixing a real number $\epsilon >0$, let $\delta_{\mathcal{F}}= \lim \sup_n \rho(\mathcal{F}_n)^{1/n}$. Then by 3.8, we can check that $\delta_{\mathcal{F}}\geq \lim_n ||M(\mathcal{F}_n)||^{1/n},$ and hence there is a positive integer $l$ such that $\dfrac{||M(\mathcal{F}_l)||}{(\delta_{\mathcal{F}}+ \epsilon)^l} r^2 ||\vec{c}|| < 1.$ We fix such $l$ and we apply the arguments of last computations to conclude, for $\mathcal{F}_{ln}=(\mathcal{F}_l)_n$, that is it true that there exists a  constant $C_1$ such that \begin{center}
$\sum_{f \in \mathcal{F}_{nl}} h_H(f(P))\leq C_1 k^{ln} n^2\dfrac{R^{n}}{(\delta_{\mathcal{F}}+ \epsilon)^{nl}}(\delta_{\mathcal{F}}+ \epsilon)^{nl}h_H(P),$  
\end{center} for $R \leq \max_{i} \{1, r^2||\vec{c} ||M(\mathcal{F}_l||\}$. Thus, there is a constant $C_2$ such that $C_1 n^2\dfrac{R^n}{(\delta_{\mathcal{F}}+ \epsilon)^{nl}} \leq C_2$ for all $n$. So we find that
\begin{center}
$\sum_{f \in \mathcal{F}_{nl}} h^+_X(f(P)) \leq C_2. k^{nl}.(\delta_{\mathcal{F}} + \epsilon)^{nl} .  h^+_X(P)$
\end{center} for all $n$, showing theorem 5.4.

\textit{Proof that theorem 5.4 implies theorem 5.1} We proved that for any $\epsilon >0$ there is a positive integer $l$ and a positive constant $C$ so that 
\begin{center}
$\sum_{f \in \mathcal{F}_{nl}} h^+_X(f(P)) \leq C. k^{nl}.(\delta_{\mathcal{F}} + \epsilon)^{nl} .  h^+_X(P),$
\end{center}for all $n,$ and $P \in X_{\mathcal{F}}(\bar{K})$. Given a integer $n$, there are $q\geq 0$ and $0<t<l$ such that $n=lq +t$. Let also $C_1$ be the constant of lemma 5.5. For $P \in X_{\mathcal{F}}(\bar{K})$, we calculate that
\begin{center}
$\sum_{f \in \mathcal{F}_n} h^+_X(f(P)) \:\:\:\:\:\:\:\:\:\:\:\:\:\:\:\:\:\:\:\:\:\:\:\:\:\:\:\:\:\:\:\:\:\:\:\:\:\newline \newline \leq C. k^{lq}.(\delta_{\mathcal{F}} + \epsilon)^{lq}. \sum_{f \in \mathcal{F}_t} h^+_X(f(P)) \newline \newline \leq C. k^{lq}.(\delta_{\mathcal{F}} + \epsilon)^{lq}.C_1^t.k^t. h^+_X(P) \:\:\:\:\:\:\:\:\:\:\newline \newline \leq CC_1^{l-1}k^n(\delta_{\mathcal{F}} + \epsilon)^{n}h^+_X(P),\:\:\:\:\:\:\:\:\:\:\:\:\:\:\:\:\:\:\:\:\:\:\:\:\:\:\:\:\:\:\:\:\:\:\:\:\:\:\:\:\:\:\:\:$
\end{center}
as we wanted to show. \newline


\section{The arising of new canonical heights}

In this final section, we show that the canonical height limit, proposed  and constructed by S. Kawaguchi in [10, theorem 1.2.1], is convergent for certain eigendivisor classes relative to algebraic equivalence, instead of linear equivalence case worked by Kawaguchi. The theorem is also an extension of theorem 5 of [12], where the eigensystem of the hypothesis has just one morphism.\newline

{\bf Theorem 6.1:} \textit{Assume that $\mathcal{F}=f_1,...,f_k:X \rightarrow X$ are morphisms, and let $D \in $Div$(X)_{\mathbb{R}}$ that satisfies the algebraic relation}
\begin{center}
$\sum^k_{i=1} f^*_iD \equiv  \beta D$\textit{ for some real number} $\beta >\sqrt{\delta_{\mathcal{F}}}k,$
\end{center}
\textit{where $\equiv$ denotes algebraic equivalence in NS$(X)_{\mathbb{R}}.$ Then }\newline

\textit{(a) For all $P \in X(\bar{K})$, the following limit converges:}
\begin{center}
$\hat{h}_{D,\mathcal{F}}(P)= \lim_{n \rightarrow \infty}\dfrac{1}{\beta^n}\sum_{ f \in \mathcal{F}_n} h_D(f(P)).$
\end{center}

\textit{(b) The canonical height in (a) satisfies }
\begin{center}
$\sum^k_{i=1} \hat{h}_{D,\mathcal{F}}(f_i(P))=\beta \hat{h}_{D,\mathcal{F}}(P)$ and $\hat{h}_{D,\mathcal{F}}(P)= h_D(P) + O(\sqrt{h^+_X(P)}). $
\end{center}

\textit{(c) If $\hat{h}_{D,\mathcal{F}}(P) \neq 0$, then $\underline{\alpha}_{\mathcal{F}}(P) \geq \beta/k.$}\newline

\textit{(d) If $\hat{h}_{D,\mathcal{F}}(P) \neq 0$ and $\beta=\delta_{\mathcal{F}}k,$ then  $\alpha_{\mathcal{F}}(P)= \delta_{\mathcal{F}}.$}\newline

\textit{(e) Assume that $D$ is ample and that $K$ is a number field. Then}
\begin{center}
$\hat{h}_{D,\mathcal{F}}(P)=0 \iff  P$ \textit{is preperiodic, i.e, has finite $\mathcal{F}$-orbit.}
\end{center}

\begin{proof}
 (a) Theorem 5.1 says that for every $\epsilon >0$ there is a constant \newline $C_1=C_1(X,h_X,\mathcal{F}, \epsilon)$ such that
\begin{center}
$\sum_{f \in \mathcal{F}_n} h^+_X(f(P)) \leq C_1.k^n. (\delta_{\mathcal{F}} + \epsilon)^n .  h^+_X(P)$ for all $n \geq 0.$
\end{center}
We are given that $\sum^k_{i=1} f^*_iD \equiv  \beta D.$ Applying lemma 5.3 with \newline  $E=\sum^k_{i=1} f^*_iD -  \beta D,$ we find a positive constant $C_2=C_2(D, \mathcal{F}, h_X)$ such that
\begin{center}
$|h_{\sum^k_{i=1} f^*_iD}(Q) - \beta h_D(Q)| \leq C_2 \sqrt{h^+_X(Q)} $ for all $Q \in X(\bar{K}).$
\end{center}
Since we assumed that the $f_i$ are morphisms, standard functoriality of Weil height states that
\begin{center}
$ h_{\sum^k_{i=1} f^*_iD} = \sum^k_{i=1} h_D \circ f_i + O(1),$
\end{center}
so the above inequality is reformulated as follows 
\begin{center}
 (**) $| \sum^k_{i=1} h_D(f_i(Q)) - \beta h_D(Q)| \leq C_3 \sqrt{h^+_X(Q)} $ for all $Q \in X(\bar{K}).$
\end{center}
For $N \geq M \geq 0$ we estimate a telescopic sum,
\begin{center}
$|\beta^{-N} \sum_{f \in \mathcal{F}_N} h_D(f(P)) - \beta^{-M} \sum_{f \in \mathcal{F}_M} h_D(f(P))| ~~~~~~~~~~~~~~~~~~ ~~~~~ ~~~~~~~~~~~~~~~\:\:\:\:\:\:\:\:\:\:\:\:\:\:\:\:\:\:\:\:\:\:\:\:\:\:\:\:\:\:\:\:\:\:\:\:\:\:\:\: \newline \newline
=|\sum^{N}_{n=M+1}\beta^{-n}[ \sum_{f \in \mathcal{F}_n} h_D(f(P))-  \beta \sum_{f \in \mathcal{F}_{n-1}} h_D(f(P))]| ~~~~~ ~~~~~~~~~~~~~~~~~~~\:\:\:\:\:\:\:\:\:\:\:\:\:\:\:\:\:\:\:\newline \newline
\leq \sum^{N}_{n=M+1}\beta^{-n}|\sum_{f \in \mathcal{F}_n} h_D(f(P))-  \beta \sum_{f \in \mathcal{F}_{n-1}} h_D(f(P))| ~~~~~~~~~~~~~~~~\:\:\:\:\:\:\:\:\:\:\:\:\:\:\:\:\:\:\:\:\:~~~~~~~~~~\newline \newline
\leq \sum^{N}_{n=M+1}\beta^{-n} [\sum_{f \in \mathcal{F}_{n-1}}|\sum^k_{i=1} h_D(f_i(f(P))) - \beta h_D(f(P))| ] ~~~~~~~~~~ ~~~~~~~~~~~~\:\:\:\:\:\:\:\:\:\:\:\:\:\:\:\:\:\:\newline \newline
\leq  \sum^{N}_{n=M+1}\beta^{-n} (\sum_{f \in \mathcal{F}_{n-1}} C_3 \sqrt{h^+_X(f(P))}) $  by (**) $~~~~~~~~~~~~~~~~~~~~~ ~~~~~~~ ~~~\:\:\:\:\:\:\:\:\:\:\:\:\:\:\:\:\:\:\:\:\:\:\:\:\:\:\:\:\:\:\:\:\:\:\:\:\:\:\:\:~~~~~\newline \newline
 \leq  \sum^{N}_{n=M+1}\beta^{-n}.k^{(n-1)/2}. C_3 . \sqrt{\sum_{f \in \mathcal{F}_{n-1}}h^+_X(f(P))}  $ by Cauchy-Schwarz ~~~~~~~~~~\:  \newline \newline
$\leq  \sum^{N}_{n=M+1}\beta^{-n}.k^{n-1}. C_3 .C. (\delta_{\mathcal{F}} + \epsilon)^{(n-1)/2}.\sqrt{h_X^+(P)}$ by Thm. 5.1 ~~~~~~~~~~~~~~~~\:\:\:\:\:\:\:\newline \newline
$\leq CC_3 \sqrt{h_X^+(P)}\sum_{n=M+1}^{\infty} [ \dfrac{k^2(\delta_{\mathcal{F}}+ \epsilon)}{\beta^{2}}]^{n/2}.$~~~~~~~~~~~~~~~~~~~~~~~~~~~~~~~~~~~~~~~~~~~~~~~~~~~~~~\:\:\:\:\:\:\:\:\:\:\:\:\:\:\:\:\:\:\:\:\:\:\:\:\:\:\:\:\:\:\:\:\:\:\:\:\:\:\:\:\:\:\:\:\:\:\:\:\:\:\:\:\:\:\:\:\:\:\:\:\:\:\:\:
\end{center} And
\begin{center}
$\sum_{n=M+1}^{\infty} [ \dfrac{k^2(\delta_{\mathcal{F}}+ \epsilon)}{\beta^{2}}]^{n/2} < \infty \iff  \dfrac{k^2(\delta_{\mathcal{F}}+ \epsilon)}{\beta^{2}} < 1.$
\end{center}
Since $\beta > \sqrt{\delta_{\mathcal{F}}k^2},$ we can choose  $0< \epsilon < \frac{\beta^{2}}{k^2} - \delta_{\mathcal{F}}$, which implies $\frac{k^2(\delta_{\mathcal{F}}+ \epsilon)}{\beta^{2}} < 1$ and the desired convergence. Also we obtain the following estimate (***):
\begin{center}
$ |\beta^{-N} \sum_{f \in \mathcal{F}_N} h_D(f(P)) - \beta^{-M} \sum_{f \in \mathcal{F}_M} h_D(f(P))| \newline \leq CC_3[ \dfrac{k^2(\delta_{\mathcal{F}}+ \epsilon)}{\beta^{2}}]^{M/2} \sqrt{h_X^+(P)}. ~~~~~~~~~~~~~~~~~~~~~~~~~~~~~~ ~~~~~~~~~~$
\end{center}
(b) The formula 
\begin{center}
$\sum^k_{i=1} \hat{h}_{D,\mathcal{F}}(f_i(P))=\beta \hat{h}_{D,\mathcal{F}}(P)$
\end{center} follows immediately from the limit defining $\hat{h}_{D,\mathcal{F}}$ in part (a). Next, letting $N \rightarrow \infty$ and setting $M=0$ in (***) gives
\begin{center}
$ |\hat{h}_{D,\mathcal{F}}(P)- h_D(P)|=  O(\sqrt{h^+_X(P)}),$
\end{center} which completes the proof of (b).

(c) We are assuming that $\hat{h}_{D,\mathcal{F}}(P) \neq 0.$ If $\hat{h}_{D,\mathcal{F}}(P) <0,$ we change $D$ to $-D,$ so we may assume $\hat{h}_{D,\mathcal{F}}(P)>0.$ Let $H \in $ Div$(X)$ be an ample divisor such that $H+D$ is also ample (this can always be arranged by replacing $H$ with $mH$ for a sufficiently large $m$). Since $H$ is ample, we may assume that the height function $h_H$ is non-negative. We compute
\begin{center}
$\sum_{f \in \mathcal{F}_n}h_{D+H}(f(P)) ~~~~~~~~~~~~~~~~~~~~~~~~~~~~~~~~~~~~~~~~~~~~~~~~~~~~~~~~~~~~~~~~~~~~~~~~~~~~~~~~~\:\:\:\:\:\:\:\:\:\:\:\:\:\:\:\:\:\:\:\:\:\:\:\:\:\:\:\:\:\:\:\:\:\:\:\:\:\:\:\:\:\:\:\:\:\:\:\:\:\:\:\:\:\:\:\:\:\:\:\:\:\:\:\:\:\:\:\:\:\:\:\:\:\:\:\:\:\:\:\:\:\:\:\:\:\:\:\:\:\:\:\:\newline \newline = \sum_{f \in \mathcal{F}_n}h_{D}(f(P)) + \sum_{f \in \mathcal{F}_n}h_{H}(f(P)) + O(k^n) ~~~~~~~~~~~~~~~~~~~~~~~~~~~~~~~~~~~~~~~~~~~~\:\:\:\:\:\:\:\:\:\:\:\:\:\:\:\:\:\:\:\:\:\:\:\:\:\:\:\:\:\:\:\:\:\:\:\:\newline \newline 
\geq  \sum_{f \in \mathcal{F}_n}h_{D}(f(P)) + O(k^n)$ since $h_H \geq0 ~~~~~~~~~~~~~~~~~~~~~~~~~~~~~~~~~~~~~~~~~~~~~~~~~~~~\:\:\:\:\:\:\:\:\:\:\:\:\:\:\:\:\:\:\:\:\:\:\:\:\:\:\:\:\:\:\:\:\:\:\:\:\:\:\:\:\:\:\:\:\:\:\:\:\:\newline \newline 
= \sum_{f \in \mathcal{F}_n}\hat{h}_{\mathcal{F},D}(f(P)) + O( \sum_{f \in \mathcal{F}_n}\sqrt{h^+_{X}(f(P))} ) $ from (b) $~~~~~~~~~~~~~~~~~~~~~~~~\:\:\:\:\:\:\:\:\:\:\:\:\:\:\:\:\:\:\:~~~~~~~~ \newline \newline
=\beta^{n}\hat{h}_{\mathcal{F},D}(P) +  O( \sum_{f \in \mathcal{F}_n}\sqrt{h^+_{X}(f(P))} ) $ from (b) $~~~~~~~~~~~~~~~~~~~~~~ ~~~~~~~~~~~~~\:\:\:\:\:\:\:\:\:\:\:\:\:\:\:\:\:\:\:\:\:\:\:\:\:\:\:\:\:\:\:\:\:\:~~\:~~~~~\newline \newline
\geq  \beta^{n}\hat{h}_{\mathcal{F},D}(P) +  O( \sqrt{\sum_{f \in \mathcal{F}_n}h^+_{X}(f(P))} ) $ since $ (x \rightarrow \sqrt{x})$ is convex $ ~~~~~~~~~~~~~~~~~~~~\:\newline \newline 
=   \beta^{n}\hat{h}_{\mathcal{F},D}(P) +  O( \sqrt{Ck^n(\delta_{\mathcal{F}} + \epsilon)^nh_X^+(P)})$ from Theorem 5.1.~~~~~~~~~~~~~~~~~~~~~~~~~\:\:\:\:\:\:\:\:\:\:\:\:\:\:\:\:\:\:\:~~~~
\end{center}
This estimate is true for every $\epsilon >0$, where $C$ depends on $\epsilon.$ Using the assumption that $\beta > \sqrt{k.\delta_{\mathcal{F}}} $ we can choose $\epsilon >0$ such that  \newline $k.(\delta_{\mathcal{F}} + \epsilon) < \beta^{2}.$ This gives
\begin{center}
$\sum_{f \in \mathcal{F}_n}h_{D+H}(f(P)) \geq  \beta^{n}\hat{h}_{\mathcal{F},D}(P) +  o(\beta^n),$
\end{center} so taking $n^{th}$-roots, using the assumption that $\hat{h}_{\mathcal{F} ,D}(P) >0,$ and letting $n \rightarrow \infty$ yields
\begin{center}
$\underline{\alpha}_{\mathcal{F}}(P)=\lim \inf_{n \rightarrow \infty} \dfrac{1}{k}\{\sum_{f \in \mathcal{F}_n}h_{D+H}(f(P) \}^{1/n} \geq \dfrac{\beta}{k}.$
\end{center}
(d) From (c) we get that $\underline{\alpha}_{\mathcal{F}}(P) \geq \dfrac{\beta}{k} =  \dfrac{\delta_{\mathcal{F}}.k}{k}=\delta_{\mathcal{F}},$ while corollary 5.2 gives $\bar{\alpha}_{\mathcal{F}}(P) \leq \delta_{\mathcal{F}}.$ Hence the limit defining $\alpha_{\mathcal{F}}(P)$ exists and is equal to $\delta_{\mathcal{F}}.$

(e) First suppose that $\# \mathcal{O}_{\mathcal{F}}(P) < +\infty.$  Since $D$ is ample and the orbit of $P$ is finite, we have that $h_D \geq 0, \hat{h}_{\mathcal{F},D}(P)  \geq 0$, and there is a constant $C>0$ such that $h_D(f(P)) \leq C $ for all $f \in \cup_{l \geq0} \mathcal{F}_l$. This gives
\begin{center}
$ |\hat{h}_{\mathcal{F},D}(P)| \leq  \lim_{n \rightarrow \infty}\dfrac{1}{\beta^n}\sum_{ f \in \mathcal{F}_n} |h_D(f(P))| \leq \lim_{n \rightarrow \infty} C.\dfrac{k^n}{\beta^n}=0$
\end{center} Since $\beta > k.$

For the other direction, suppose that $\hat{h}_{\mathcal{F},D}(P)=0.$ Then for any $n \geq 0$ and $g \in \mathcal{F}_n,$ we apply part (b) to obtain
\begin{center}
$0=\beta^n\hat{h}_{\mathcal{F},D}(P)=\sum_{f \in \mathcal{F}_n}\hat{h}_{\mathcal{F},D}(f(P)) \geq \hat{h}_{\mathcal{F},D}(g(P))  \newline \geq h_D(g(P)) - c\sqrt{h_D(g(P))}. ~~~~~~~~~~~~~~~~~~~~~~~~~~~~~~~~~~~~~~~~~$
\end{center} This gives $h_D(g(P)) \leq c^2,$ where $c$ does not depend on $P$ or $n.$ This shows that $\mathcal{O}_{\mathcal{F}}(P)$ is a set of bounded height with respect to an ample height. Since $\mathcal{O}_{\mathcal{F}}(P)$ is contained in $X(K(P))$ and since we have assumed that $K$ is a number field, we conclude that $\mathcal{O}_{\mathcal{F}}(P)$ is finite.
\end{proof} 

{\bf Remark 6.2:} In the same way as pointed in remark 29 of [12], when $f_1,...,f_k$ are morphisms, there is always one divisor class $D \in$ NS$(X)_{\mathbb{R}}$ such that $\sum^k_{i=1} f^*_iD \equiv  \beta D$, where $\beta$ is the spectral radius of the linear map $\sum_{i \leq k} A(f_i)$ on NS$(X)_{\mathbb{R}}$. It would remain to check whether it satisfies $\beta >k.\sqrt{\delta_{\mathcal{F}}}$. This works for the non-trivial example 3.3, 3.4, 3.5 and 3.6, where the above height will coincide with the height constructed by Kawaguchi and Silverman.


\begin{thebibliography}{9}                                                                                                %


\bibitem{}  A. Baragar,\ \textit{Canonical vector heights on algebraic K3 surfaces with Picard number two}, Canad. Math. Bull. 46 (2003), 495-508 .

\bibitem{}  A. Baragar,\ \textit{Rational points on K3 surfaces in $\mathbb{P}^1 \times\mathbb{P}^1 \times \mathbb{P}^1$ },Math. Ann. 305 (1996),541-558 .


\bibitem{}  E. Bombieri and W. Gubler,\ \textit{ Heights in Diophantine Geometry}, Number 4 in New Mathematical Monographs. Cambridge University Press, Cambridge, 2006.

\bibitem{} G. Call and H. Silverman, \ \textit{Canonical heights on varieties with morphisms},
Compositio Math. 89 (1993), 163-205. 


\bibitem{}  W. Fulton, \ \textit{Intersection theory, 2nded.,Ergeb.Math.Grenzgeb.(3)2}, Springer-Verlag, Berlin 1998.

\bibitem{}  V.Guedj, \ \textit{Ergodic properties of rational mappings with large topological degree}, Ann.of Math.(2) 161 (2005), no. 3, 1589-1607.

\bibitem{}  R. Hartshorne, \ \textit{Algebraic Geometry}, Springer-Verlag, New
York, 1977.

\bibitem{}  M. Hindry and J. H. Silverman, \ \textit{Diophantine Geometry: An Introduction}, volume 201 of Graduate Texts in Mathematics. Springer-Verlag, New
York, .

\bibitem{}  S. Kawaguchi, \ \textit{Canonical height functions for affine plane automorphisms}, Math. Ann., 335(2):285-310, 2006.

\bibitem{}  S. Kawaguchi, \textit{Canonical heights, invariant currents, and dynamical eigensystems of morphisms for line bundles,} preprint, math.NT/0405006. (The revised version is accepted in J. Reine Angew. Math.)

\bibitem{}  S. Kawaguchi, \ \textit{Projective surface automorphisms of positive topological entropy from an arithmetic viewpoint},  available at arXiv:math/0510634v2 [math.AG].

\bibitem{}  S. Kawaguchi and J. H. Silverman, \ \textit{On the dynamical and arithmetic
degrees of rational self-maps of algebraic varieties},2012. arXiv:1208.0815.

\bibitem{} S. Kawaguchi and J. H. Silverman, \ \textit{Dynamical canonical heights for Jordan blocks, arithmetic degrees of orbits, and nef canonical heights on abelian varieties}, preprint 2013, http://arxiv.org/abs/1301.4964




\bibitem{}  S. Lang, \ \textit{Fundamentals of Diophantine Geometry},  New York, (1983).


\bibitem{} Y. Matsuzawa, \ \textit{On upper bounds of arithmetic degrees}, preprint 2016, http://arxiv.org/abs/1606.00598



\bibitem{}  J.H. Silverman, \ \textit{Examples of dynamical degree equals arithmetic degree}, available at arXiv:1212.3015v1[math.NT].

\bibitem{}  J.H. Silverman, \ \textit{Dynamical Degrees, Arithmetic Entropy, and Canonical Heights for Dominant Rational Self-Maps of Projective Space}, available at arXiv:1111.5664[math.NT].

\bibitem{}  J.H. Silverman, \ \textit{Heights and the specialization maps for families of abelian varieties},  J. Reine Angew. Math. 342(1983) 197-211.

\bibitem{} J.H. Silverman, \ \textit{Rational points on K3 surfaces: a new canonical height}, Invent. Math. 105 (1991), 347-373.

\bibitem{}  J.H.Silverman, \ \textit{The Arithmetic of Dynamical Systems}, volume 241 of Graduate Texts in Mathematics. Springer, New York, 2007.

\bibitem{} J.H.Silverman, \ \textit{The Arithmetic of Elliptic Curves}, 2nd ed., Grad. Texts in Math. 106, Springer-Verlag, Dordrecht 2009.





\end{thebibliography}
\end{document}